\documentclass[12pt]{article}
\usepackage[usenames]{color}
\usepackage{colortbl}
\usepackage{tikz-cd}
\usepackage{xcolor}
\usepackage[bbgreekl]{mathbbol}

\usepackage{amscd,amsmath, amssymb, fancyhdr, mathbbol, dutchcal, url, euscript, amssymb, mathrsfs}
\usepackage[backref=page]{hyperref}
\usepackage{mathtools}
\renewcommand*{\backrefalt}[4]{%
	\ifcase #1 (Not cited.)%
	\or        (Cited on page~#2.)%
	\else      (Cited on pages~#2.)%
	\fi}

\hypersetup{
	colorlinks   = true,
	citecolor    = magenta,
	linkcolor    =blue,
	urlcolor     =magenta	
}

\numberwithin{equation}{section}
\newcommand{\version}{version 1.0,\ \ October, 2023}

\def\eqref#1{(\ref{#1})}

\newcommand{\g}{{\mathfrak g}}

\newcommand{\arrow}{{\:\longrightarrow\:}}
\newcommand{\Z}{{\Bbb Z}}
\def\C{{\Bbb C}}
\newcommand{\R}{{\Bbb R}}
\newcommand{\Q}{{\Bbb Q}}
\renewcommand{\H}{{\Bbb H}}

\def\1{\sqrt{-1}\:}

\newcommand{\cntrct}                
{\hspace{2pt}\raisebox{1pt}{\text{$\lrcorner$}}\hspace{2pt}}

\renewcommand{\tilde}{\widetilde}
\renewcommand{\bar}{\overline}
\renewcommand{\phi}{\varphi}
\renewcommand{\epsilon}{\varepsilon}
\renewcommand{\geq}{\geqslant}

\newcommand{\Grad}{\nabla\!\!\!\!\nabla}

\newcommand{\End}{\operatorname{End}}

\newcommand{\Id}{\operatorname{Id}}

\newcommand{\Hom}{\operatorname{Hom}}

\newcommand{\Lie}{\operatorname{Lie}}

\newcommand{\Spec}{\operatorname{Spec}}

\newcommand{\Mat}{\operatorname{Mat}}

\renewcommand{\Im}{\operatorname{Im}}
\newcommand{\Aff}{\operatorname{Aff}}
\newcommand{\GL}{\operatorname{GL}}

\newcounter{Mycounter}[section]
\newcounter{lemma}[section]
\newcounter{claim}[section]
\renewcommand{\theclaim}{{Claim \thesection.\arabic{claim}}}
\newcommand{\claim}{%
    \setcounter{claim}{\value{Mycounter}}
    \refstepcounter{claim}
    \stepcounter{Mycounter}
    {\noindent \bf \theclaim:\ }}

\newcounter{sublemma}[section]
\newcounter{corollary}[section]
\renewcommand{\thecorollary}{{Corollary \thesection.\arabic{corollary}}}
\newcommand{\corollary}{%
    \setcounter{corollary}{\value{Mycounter}}
    \refstepcounter{corollary}
    \stepcounter{Mycounter}
    {\noindent \bf \thecorollary:\ }}

\newcounter{theorem}[section]
\renewcommand{\thetheorem}{{Theorem \thesection.\arabic{theorem}}}
\newcommand{\theorem}{%
    \setcounter{theorem}{\value{Mycounter}}
    \refstepcounter{theorem}
    \stepcounter{Mycounter}
    {\noindent \bf \thetheorem:\ }}

\newcounter{conjecture}[section]
\renewcommand{\theconjecture}{{Conjecture \thesection.\arabic{conjecture}}}
\newcommand{\conjecture}{%
    \setcounter{conjecture}{\value{Mycounter}}
    \refstepcounter{conjecture}
    \stepcounter{Mycounter}
    {\noindent \bf \theconjecture:\ }}

\newcounter{proposition}[section]
\newcounter{definition}[section]
\renewcommand{\thedefinition}
      {{Definition~\thesection.\arabic{definition}}}
\newcommand{\definition}{%
    \setcounter{definition}{\value{Mycounter}}
    \refstepcounter{definition}
    \stepcounter{Mycounter}
    {\noindent \bf \thedefinition:\ }}

\newcounter{example}[section]
\renewcommand{\theexample}{{Example \thesection.\arabic{example}}}
\newcommand{\example}{%
    \setcounter{example}{\value{Mycounter}}
    \refstepcounter{example}
    \stepcounter{Mycounter}
    {\noindent \bf \theexample:\ }}

\newcounter{remark}[section]
\renewcommand{\theremark}{{Remark \thesection.\arabic{remark}}}
\newcommand{\remark}{%
    \setcounter{remark}{\value{Mycounter}}
    \refstepcounter{remark}
    \stepcounter{Mycounter}
    {\noindent \bf \theremark:\ }}

\newcounter{properties}[section]
\renewcommand{\theproperties}{{Properties \thesection.\arabic{properties}}}
\newcommand{\properties}{%
    \setcounter{properties}{\value{Mycounter}}
    \refstepcounter{properties}
    \stepcounter{Mycounter}
    {\noindent \bf \theproperties:\ }}

\newcounter{problem}[section]
\newcounter{question}[section]
\newcommand{\proof}{\noindent{\bf Proof:\ }}

\makeatletter

\@addtoreset{equation}{section} \@addtoreset{footnote}{section}
\makeatother

\def\blacksquare{\hbox{\vrule width 5pt height 5pt depth 0pt}}
\def\endproof{\blacksquare}

\title{Flat hypercomplex nilmanifolds are $\H$-solvable}
\author{Yulia Gorginyan}
\begin{document}
\maketitle
\begin{abstract} 
We say that a hypercomplex nilpotent Lie algebra is {\bf $\H$-solvable} if there exists a sequence of $\H$-invariant subalgebras 
\begin{equation*}
 \g_1^{\H}\supset\g_2^{\H}\supset\cdots\supset\g_{k-1}^{\H}\supset\g_k^{\H}=0,
\end{equation*}
such that $[\g_i^{\H},\g_i^{\H}]\subset\g^{\H}_{i+1}.$ Let $N=\Gamma\backslash G$ be a hypercomplex nilmanifold with flat Obata connection and $\g=\Lie(G)$. We prove that the Lie algebra $\g$ is $\H$-solvable.
\end{abstract}
\tableofcontents


\section{Introduction}
\subsection{Affine manifolds}

A manifold $X$ together with a torsion-free flat connection $\nabla$ is called {\bf an affine manifold}.
Equivalently, an affine manifold $X$ is a manifold with an atlas such that all translation maps between charts are in $\Aff(\R^n)$ (\cite{FGH}, \cite{Sh}). 

\medskip
Recall that a group of affine transformations $\Aff(\R^n)$ is a semidirect product $\GL_n(\R)\ltimes\R^n$. If $f\in\Aff(\R^n)$ then for any $x\in\R^n$ there exist a linear part $A\in\GL(\R^n)$ and a translation part $t\in\R^n$ such that $f(x)=Ax+t$. {\bf The linearization map} is the natural homomorphism $l: \Aff(\R^n)\arrow\GL(\R^n)$ given by the formula $l(f)=A$. 

\medskip
An affine manifold $X$ is called {\bf complete} if its universal cover is affine equivalent to $\R^n$, i.e. $X=\R^n/\Gamma$, where $\Gamma\subset\Aff(\R^n)$ is a discrete subgroup of the affine transformation group. 


\medskip
There are three famous open questions related to affine geometry. In the historical context, the first was the conjecture of Chern:

\medskip
\conjecture (Chern) The Euler characteristic of a compact affine manifold vanishes \cite{Gol}.

\medskip
Kostant and Sullivan proved Chern's conjecture for the compact complete affine manifolds \cite{KS}. Also, it was proven by Bruno Klingler \cite{Klin} in the case when a compact affine manifold admits a parallel volume form.
As far as we know, the general case of Chern’s conjecture remains open. 

\medskip
The conjecture of Markus links the existence of a parallel volume form with the completeness of a compact affine manifold.

\medskip
\conjecture (Markus, 1962)\cite{Mar}
A compact affine manifold $X$ admits parallel volume form if and only if the manifold $X$ is complete.



\medskip
Recall that a subgroup $\Gamma$ of an affine group $\Aff(\R^n)$ is called {\bf crystallographic} \cite{A} if its action on $\R^n$ is properly discontinuous, free, and cocompact. 

\medskip
Recall that a group $G$ is said to have a certain property $P$ {\bf virtually}
if $G$ contains a subgroup $H$ of finite index which has the property $P$.

\medskip
\begin{theorem}(Bieberbach, 1911) \cite{Bieb} Every discrete subgroup $\Gamma$ of an isometry group ${\rm Isom}(\R^n)$ is virtually abelian. Every crystallographic subgroup of ${\rm Isom}(\R^n)$ is virtually a translation group. For a given $n$ there exists only a finite number of such groups $\Gamma\subset{\rm Isom}(\R^n)$ up to a conjugation.
\end{theorem}

\medskip
One of the ways to generalize the theorem of Bieberbach is
to consider an affine transformation group $\Aff(\R^n)$ instead of ${\rm Isom}(\R^n)$. 

\medskip
\conjecture (Auslander, 1964)\cite{Aus} Every crystallographic subgroup of an affine group is virtually solvable, i.e. contains a solvable subgroup of finite index.

\medskip
Auslander’s conjecture was proven by D. Fried and W. Goldman \cite{FG} in the case $n=3$ and the result was refined up to dimension $n=6$ by H. Abels, G.A. Margulis, and G.A. Soifer \cite{AMS}.
Also, Auslander’s conjecture is true in a case when the monodromy group preserves a metric of signature $(1,n)$ \cite{GK} and $(2,n)$ \cite{AMS2}.




\subsection{Hypercomplex affine nilmanifolds}

Recall that {\bf a nilmanifold} $N$ is a compact quotient of a connected simply connected nilpotent Lie group $G$ by a lattice subgroup $\Gamma$, which acts on the group $G$ from the left. We denote the quotient as $N=\Gamma\backslash G$.
 It is called {\bf an affine nilmanifold} if $G$ has a left-invariant affine structure, i.e. the flat torsion-free connection $\nabla$ such that $L_g^*(\nabla)=\nabla$, where $L_g: G\arrow G$ is the left-translation.

Let $\pi:\tilde X\arrow X$ be the universal cover of a manifold $X$ and $\pi_1(X)$ the fundamental group. An affine immersion $D:\tilde{X}\arrow\R^n$ is called {\bf the developing map}. By a well-known theorem (\cite{OV}, \cite{Gol1}) $D$ exists and is uniquely defined up to an affine automorphism.

\medskip
\begin{definition} {\bf The affine holonomy representation} is a
unique homomorphism $h:\pi_1(X)\arrow\Aff(\R^m)$ which satisfies $D\circ \gamma=h(\gamma)\circ D$ for every $\gamma\in\pi_1(X)$. 
{\bf The affine holonomy group} $\mathcal{H}:=\Im(h)$ is the image of the homomorphism $h$ in $\Aff(\R^n)$.
{\bf The linear holonomy group} is the image of the affine holonomy group under the linearization map, $\mathcal{L}:=l(\mathcal{H})\subset\GL(\R^n)$.
\end{definition}

\medskip
\begin{definition}
A closed Lie subgroup $G\subset\GL(V)$ is called {\bf a linear algebraic group} if $G$ is given by a system of polynomial equations.
\end{definition}

\medskip

\begin{definition}\label{unipotent}
A representation of a linear algebraic group is called {\bf unipotent} if, in a certain basis, its image is the group of strictly upper-triangular matrices.
\end{definition}



\medskip 
Let $X$ be a compact affine manifold whose linear holonomy representation is unipotent. Then
$X$ admits a parallel volume form. The converse was proven by Goldman, Fried, and Hirsch, \cite[Theorem A]{FGH}: 

\medskip
\begin{theorem}\label{FGH'}
Let $X$ be a compact affine manifold with a parallel volume form. Assume that its affine holonomy group is nilpotent. Then the linear holonomy representation is unipotent.
\end{theorem}

\medskip
We will use the following reformulation of \ref{FGH'}:

\medskip
\begin{theorem}\label{FGH}
  {Let $X$ be a compact affine manifold with a parallel volume form. Assume that its fundamental group is nilpotent. Then the monodromy representation is unipotent.}
\end{theorem}

\medskip
We are interested in affine nilmanifolds which also possess a hypercomplex structure. First, recall a definition of a complex nilmanifold.

\medskip
\begin{definition}\label{comnil}
Let $G$ be a Lie group equipped with a left-invariant complex structure. {\bf A complex nilmanifold} is a pair $(N=\Gamma\backslash G, I)$, where $N$ is a nilmanifold and the complex structure $I$ is obtained from the corresponding left-invariant complex structure on $G$.
\end{definition}

\medskip
An almost hypercomplex manifold $X$ is a smooth manifold equipped with three endomorphisms $I, J$ and $K$ of the tangent bundle satisfying the quaternionic relations $I^2=J^2=K^2=-\Id$ and $IJ=K$. When almost complex structures are integrable, the quadruple $(X, I, J, K)$ is called {\bf a hypercomplex manifold.}

\medskip
\begin{definition}
Let $G$ be a nilpotent Lie group with a left-invariant hypercomplex structure $I, J, K$ and
$\Gamma\subset G$ a cocompact lattice. Then $(N=\Gamma\backslash G, I, J, K)$ is called {\bf a hypercomplex nilmanifold}.
\end{definition}

\medskip
M. Obata showed \cite{Ob} that on a hypercomplex manifold $M$, there exists a unique torsion-free connection $\nabla$ preserving the complex structures: $\nabla I=\nabla J=\nabla K=0$. It is called {\bf the Obata connection}.

It could be written in the following form
\begin{equation}\label{Obata}
    \nabla_XY=\frac{1}{2}([X,Y]+I[IX,Y]-J[X,JY]+K[IX,JY]),
\end{equation}
for any $X,Y\in TM$ \cite{Sol}.

\medskip
The existence of a parallel volume form on a hypercomplex nilmanifold is guaranteed by the following theorem: 

\medskip
\begin{theorem}\cite[Theorem 3.2]{BDV}\label{V}
Let $N=\Gamma\backslash G$ be a hypercomplex nilmanifold, $n=\dim_{\C}G$. Then $G$ admits a left-invariant, non-zero, holomorphic section $\Omega$ of the canonical bundle $\Lambda^{n,0}G$. Moreover, $\nabla\Omega=0$, where $\nabla$ is {\bf the Obata connection}.
\end{theorem}

\subsection{$\H$-solvable Lie algebras and algebraic holonomy}

An operator $I$ on a real Lie algebra $\g$ is called {\bf a complex structure operator} if $I^2=-\Id$ and the $\sqrt{-1}$-eigenspace $\g^{1,0}$ is a Lie subalgebra in the complexification $\g_{\C}=\g\otimes\C$.

\medskip
\remark
Take an endomorphism $I\in\End\g$, $I^2=-\Id$, and extend it to the left-invariant almost complex structure on the Lie group $G$. Then $I$ is complex if and only if $[\g^{1,0},\g^{1,0}]\subset\g^{1,0}$. In other words, this definition is compatible with \ref{comnil}.

\medskip
{\bf A hypercomplex structure} on a Lie algebra $\g$ is a triple of complex structure operators $I,J$ and $K$ on $\g$ satisfying the quaternionic relations.

\medskip
\begin{definition}
Let $\g$ be a nilpotent hypercomplex Lie algebra.
Define inductively $\H$-invariant Lie subalgebras:
$$\g_i^{\H}:=\H[\g_{i-1}^{\H},\g_{i-1}^{\H}],$$
where $\g_0^{\H}=\g$ and $$\g_1^{\H}:=\H[\g,\g]=[\g,\g]+I[\g,\g]+J[\g,\g]+K[\g,\g].$$
It is natural to study if $\g_1^{\H}$ is a proper subalgebra in $\g$.

\medskip
A hypercomplex nilpotent Lie algebra $\g$ is called {\bf $\H$-solvable} if the following sequence terminates for some $k\in\Z_{>0}$:
\begin{equation}\label{Hsol}
 \g_1^{\H}\supset\g_2^{\H}\supset\cdots\supset\g_{k-1}^{\H}\supset\g_k^{\H}=0.
\end{equation}
\end{definition}

\medskip
\conjecture\label{conj} Let $(\g, I, J, K)$ be a nilpotent hypercomplex Lie algebra. Then it is $\H$-solvable.

\medskip
Let $(N=\Gamma\backslash G, I, J, K)$ be a hypercomplex nilmanifold.
Consider a complex nilmanifold $(N, L)$ with a general complex structure $L=aI+bJ+cK$, where $(a,b,c)\in S^2$, obtained from the hypercomplex structure. It is natural to try to describe complex submanifolds in $(N, L)$. We solved this problem in the case of complex curves under the additional assumption on the Lie algebra $\g=\Lie(G)$  \cite{Gor}. Precisely, if the corresponding Lie algebra $\g$ is $\H$-solvable, then there are no complex curves in a complex nilmanifold $(N, L)$ for the general complex structure $L$. 

In this work, we show that in the case when a hypercomplex nilmanifold $\Gamma\backslash G$ admits the flat Obata connection, then $\g=\Lie(G)$ is $\H$-solvable. The main argument we use relies on the notion of {\it an algebraic holonomy group}, which we define below \eqref{alghol}.

\medskip
\example
It was shown in \cite{DF} that all 8-dimensional hypercomplex nilpotent Lie algebras are Obata-flat. However, in the same paper, I. Dotti and A. Fino presented an example of the 3-step nilpotent hypercomplex Lie algebra of dimension 12 which has non-zero curvature of the Obata connection.

\medskip
Let $\g$ be a nilpotent Lie algebra over the field $\R$. {\bf A rational structure} in $\g$ is a subalgebra $\g_{\Q}\subset\g$ defined over the rational numbers such that  $\g_{\Q}\otimes_{\Q}\R=\g$. Assume that a Lie algebra $\g$ admits a rational structure $\g_{\Q}$. By \cite{CG} it happens if and only if there exists a nilmanifold $\Gamma\backslash G$ such that $\g_{\Q}:=\rm{span}_{\Q}\langle\log\Gamma\rangle$. 


In the case when a Lie group $G$ admits a left-invariant hypercomplex structure with the flat Obata connection, we prove that $\g_i^{\H}$ is a proper subalgebra of $\g_{i-1}$ using the following approach.

\medskip
{\bf A connection} on a Lie algebra $\g$ is an $\R$-linear map 
\begin{equation*}
    \nabla:\g\arrow\g^*\otimes\g
\end{equation*}
such that $X\mapsto\nabla_X\in\End(\g)$. Note that a connection on a Lie algebra $\g$ is the same as a left-invariant connection on a Lie group $G$ \eqref{leftconn}. This notion can be generalized to arbitrary $G$-equivariant vector bundle on a Lie group $G$, giving the algebraic version of an equivariant connection (see \eqref{3.1}).

\medskip
{\bf The curvature tensor  $R$} of the connection $\nabla$ defined in the following way: $$R(X, Y)=[\nabla_X,\nabla_Y]-\nabla_{[X, Y]},$$ where $X, Y\in\g$. The connection  $\nabla$ is called {\bf flat} if $R=0$. Notice that left invariant flat connections on $G$ are equivalent to the representations of a Lie algebra on itself considered as a vector space \eqref{repflat}. 



\medskip
In the end, we prove the following theorem:

\medskip
\begin{theorem}
Let $(N, I, J, K)$ be a hypercomplex nilmanifold with the flat Obata connection. Then the corresponding Lie algebra is $\H$-solvable.
\end{theorem}


\section{Preliminaries: Nilpotent Lie groups and algebras}

Let $G$ be a Lie group. Define inductively the descending chain of normal subgroups of $G$
\begin{equation}\label{NG}
    G=G_0\supset G_1\supset G_2 \supset\cdots\supset G_k \supset\cdots,
\end{equation}
where $G_j:=[G_{j-1},G]$ is the subgroup generated by the elements of the form $xyx^{-1}y^{-1},\, x\in G_{j-1},\,y\in G$.

\medskip
\begin{definition}
A Lie group $G$ is called {\bf a nilpotent Lie group} if \eqref{NG} terminates for some $k\in\Z_{>0}$, i.e. $G_k=\{e\}$.
\end{definition}

\medskip
 {\bf The descending series} of a real Lie algebra $\g$ is the chain of ideals defined as follows:
 \begin{equation*}
   \g_0\supset\g_1\supset\dots\supset[\g_k,\g]\supset\dots,  
 \end{equation*}
 $\g_0=\g$ and $\g_k=[\g_{k-1},\g]$. 
 
\medskip
\begin{definition}\label{nilpotent}
A Lie algebra $\g$ is called {\bf nilpotent} if $\g_s=0$ for some $s\in\Z_{>0}$.
\end{definition}

\medskip
Let $\g$ be a real Lie algebra and $\g^*$ its dual space. For any $\alpha\in\g^*$ and $\xi,\theta\in\g$ the Chevalley–Eilenberg  differential $d:\g^*\arrow\Lambda^2\g^*$ is defined as follows 
\begin{equation}\label{ceform}
    d\alpha(\xi,\theta)=-\alpha([\xi,\theta]).
\end{equation}
It extends to a finite-dimensional complex 
\begin{equation}\label{CH}
    0\arrow\g^*\arrow\Lambda^2\g^*\arrow\cdots\arrow\Lambda^{2n}\g^*\arrow 0
\end{equation}
by the Leibniz rule: $d(\alpha\wedge\beta)=d\alpha\wedge\beta+(-1)^{\tilde{\alpha}}\alpha\wedge d\beta$, where $\alpha,\beta\in\g^*$. The identity $d^2=0$ follows from the Jacobi identity \cite{CE}. 

\medskip
Notice that the kernel of a closed $1$-form is an ideal in the Lie algebra $\g$.
According to the definition of the Chevalley–Eilenberg differential, any closed 1-form $\alpha\in\g^*$ vanishes on the commutator ideal  $[\g,\g]\subseteq\ker\alpha$. 

\medskip
\remark\label{intersectonofkernels}
The intersection of all kernels of closed 1-forms $\alpha\in\g^*$ $$\Sigma\,\,:=\bigcap_{\alpha\in\g^*, d\alpha=0}{\ker\alpha}$$
also forms an ideal in the Lie algebra $\g$ which coincides with the commutator ideal $[\g,\g]$.

\medskip
Recall that {\bf a distribution} on a smooth manifold $N$ is a sub-bundle $\Sigma\subset TN$. The distribution is called {\bf involutive} if it is closed under the Lie bracket. {\bf A leaf} of the distribution $\Sigma$ is a maximal connected, immersed submanifold $L\subset N$ such that $L$ tangent to $\Sigma$ at each point.  If\, $\Sigma$ is involutive, then the set of all its leaves is called {\bf a foliation}.

\section{Algebraic holonomy group}

\subsection{Left equivariant vector bundles}\label{3.1}
Recall that a map $\varphi:X\arrow Y$ of two manifolds $X$ and $Y$ with an action of a group $G$ is called {\bf $G$-equivariant} if $\varphi(g\cdot x)=g\cdot\varphi(x)$ for any $g\in G$, $x\in X$. A vector bundle $\pi:\mathbb{B}\arrow X$ on a manifold $X$ is called {\bf a $G$-equivariant vector bundle} if its total space is equipped with a $G$-action such that the projection $\pi$ is a $G$-equivariant map and the action of $G$ is linear on the fibers, i.e. $\mathbb{B}_{x}\arrow\mathbb{B}_{gx}$ is a linear map for all $g\in G$.

Let $G$ be a Lie group, $L_g:G\arrow G$ the left translation ($h\mapsto g\cdot h$), and $B$ a finite-dimensional vector space with a basis $\{b_1,\dots,b_n\}$. 

\medskip
\begin{definition}\label{defalgconn} Let $\g=\Lie(G)$ be a Lie algebra and $\g^*=\Hom(\g,\R)$ its dual.
    {\bf An algebraic connection} is an $\R$-linear map $$\nabla:B\arrow\g^*\otimes B.$$ We use the following notation
\begin{equation}\label{nabla}
    \nabla(b)=\sum_i\theta_i\otimes A_i(b)
\end{equation}
{ where $\theta_i\in\g^*$ and $A_i\in\End(B)$, $b\in B$.}\\
\end{definition}

 Let $d:\g^*\arrow\Lambda^2\g^*$ be the Chevalley–Eilenberg differential.  We define {\bf the curvature} of the connection $\nabla$ by the following formula
\begin{equation}\label{curv1}
\Theta_{\nabla}:=d\omega+\omega\wedge\omega,
\end{equation}
{where $\omega:=\sum_i\theta_i\otimes A_i$, $d\omega=\sum_id\theta_i\otimes A_i$ and $\omega\wedge\omega=\sum_{i<j}\theta_i\wedge\theta_{j}\otimes[A_i,A_j]$.}\\


Consider a vector bundle $\pi:\mathbb{B}\arrow G$ on a Lie group $G$. From now we denote the fiber at the identity $e\in G$ by $\mathbb{B}_e=B$. 

\medskip
\begin{definition}
    {\bf A left equivariant vector bundle} $\pi:\mathbb{B}\arrow G$ over a Lie group $G$ is a $G$-equivariant vector bundle $\pi$ on a Lie group $G$ with an action given by the left translations. {\bf A left invariant section} $s$ is determined by its value $s_e$ at $e\in G$. This gives a map
    \begin{equation}\label{leftsection}
  \mathbb{L}:\mathbb{B}_e\arrow H^0(G,\mathbb{B}),\quad  s=\mathbb{L}(s_e),
    \end{equation}
    $$\mathbb{L}(s_e)(g)=s(g):=(L_g)_*(s_e).$$
  identified the fiber $\mathbb{B}_e$ with the space $H^0(G,\mathbb{B})$ of all $G$-invariant sections of $\mathbb{B}$.
\end{definition}

\medskip
\remark
All tensor powers of a $G$-equivariant vector bundle are $G$-equivariant.

\medskip
\claim
The category of left equivariant vector bundles on a Lie group $G$ is equivalent to the category of vector spaces.

\medskip
\proof
The group $G$ acts on itself by the left multiplication freely and transitively. Any left equivariant bundle $\mathbb{B}$ is trivial as a vector bundle and there is a basis of left invariant sections $\{s_1,\dots,s_n\}$, such that $s_i=\mathbb{L}(b_i)$. For each left equivariant vector bundle the vector space of its left invariant sections gives a functor to the category of vector spaces.
Conversely, given a vector space, consider it as a space of 
sections evaluated at the identity $e\in G$, and then extend them via all left translations \eqref{leftsection} to obtain a left equivariant vector bundle.  \endproof

Now we are able to give a definition of an invariant connection on the bundle $\mathbb{B}$. We use an algebraic connection \eqref{nabla} on the fiber $B=\mathbb{B}_e$.

\medskip
\begin{definition} Let $\pi:\mathbb{B}\arrow X$ be a $G$-equivariant vector bundle over a manifold $X$.
   A connection $\Grad$ on $\mathbb{B}$ is called {\bf equivariant} if it defines a left invariant differential operator $$\Grad:\mathbb{B}\arrow\Lambda^1(X)\otimes\mathbb{B}.$$ 
\end{definition}

\medskip
\claim\label{leftconn}
 A connection $\Grad:\mathbb{B}\arrow\Lambda^1(G)\otimes\mathbb{B}$ on a left equivariant vector bundle $\mathbb{B}$ over a Lie group $G$ is {equivariant} if it satisfies
 \begin{equation}\label{l-e connection}
 \Grad(\mathbb{L}(s_e))=\mathbb{L}((\Grad s)_e),
 \end{equation}
 where $\Grad s$ is understood as a section of the left-equivariant bundle $\Lambda^1(G)\otimes\mathbb{B}$, $\mathbb{L}$ is the left-translation map defined in \eqref{leftsection}, and $(\Grad s)_e$ means the value of the section in $e\in G$.


\begin{definition} Let $\theta\in\g^*$ and $b\in B$. Define a map $$d_{\nabla}:\g^*\otimes_{\R} B\arrow\Lambda^2\g^*\otimes_{\R} B$$ by the following formula:
\begin{equation}
   d_{\nabla}(\theta\otimes b)=d\theta\otimes b-\theta\wedge\nabla b.
\end{equation}
\end{definition}



This defines {\bf the curvature of an algebraic connection} $$\Theta=d^2_{\nabla}: B\arrow \Lambda^2\g^*\otimes B.$$
As in \eqref{curv1} one can also write $\Theta=d\omega+\omega\wedge\omega$.

\medskip
\claim\label{curvature}
    The curvature of a left equivariant connection $\Grad$ on a G-equivariant bundle $\pi:\mathbb{B}\arrow X$
   is an $\R$-linear map $\Theta_{\Grad}:\mathbb{B}\arrow\Lambda^2(X)\otimes\mathbb{B}$
 given by $\Theta_{\Grad}(s):=\mathbb{L}(d^2_{\nabla}(s_e))$.   
 \endproof  

\medskip
   We say that bundle $\mathbb{B}$ is {\bf flat} if $\Theta_{{\Grad}}=0$.

\medskip
\claim\label{repflat} Let $\pi:\mathbb{B}\arrow G$ be a left equivariant vector bundle over a Lie group. Then the bundle is flat if and only if $\nabla_X:\g\arrow\End(B)$ is a Lie algebra representation.

\medskip
\proof Assume that $\nabla:\g\arrow B^*\otimes B\cong\End(B)$ is a representation and the map is given as follows $X\mapsto\nabla_X$. Then for any $X,Y\in\g$ we have$[\nabla_X,\nabla_Y]-\nabla_{[X,Y]}=0$, hence by \ref{curvature} the curvature vanishes. Conversely, if $\Theta_{\Grad}=0$, then $\Theta_{\nabla}=[\nabla_X,\nabla_Y]-\nabla_{[X,Y]}=0$, hence $\nabla_X$ is a representation.
\endproof

\medskip 
Consider a left-equivariant vector bundle $\pi:\mathbb{B}\arrow G$ over a Lie group  $G$.
Let $\Gamma$ be a discrete subgroup of a Lie group $G$ which acts on $G$ from the left and $q: G\arrow\Gamma\backslash G$ the quotient map. Denote by $\mathbb{B}_{\Gamma}$ the induced vector bundle $\pi_{\Gamma}:\mathbb{B}_{\Gamma}\arrow \Gamma\backslash G$ on a manifold $\Gamma\backslash G$ obtained by taking the fiberwise quotients by $\Gamma$ such that the following diagram is commutative
\begin{center}
    \begin{tikzcd}
  \mathbb{B} \arrow[r, "q_{\Gamma}"] \arrow[d, "\pi"]
    & \mathbb{B}_{\Gamma} \arrow[d, "\pi_{\Gamma}" ] \\
  G \arrow[r, "q" ]
& || {\Gamma\backslash G}.
\end{tikzcd}  
\end{center}

Given a vector bundle $(\mathbb{B},\Grad)$ with an equivariant connection $\Grad$, we denote by $\Grad^{\Gamma}$ the induced connection on a bundle $\mathbb{B}_{\Gamma}$.
 


\subsection{Algebraic holonomy group}
Recall that there is a bijection between the set of isomorphism classes of flat vector bundles over $X$ and the set of conjugacy classes of homomorphisms $\varphi:\pi_1(X)\arrow\GL(\R^n)$ (It is called  ``the Riemann-Hilbert correspondence'', see e.g. \cite{OV}).

\medskip
Let $p: E\arrow M$ be a vector bundle over a manifold $M$ and $\delta:[0,1]\arrow M$ a smooth path in $M$. An equation $\nabla_{\dot\delta(t)}v=0$ defines {\bf a parallel transport} of a section $v$ of $E_{|_{\delta(t)}}$ along the curve $\delta$.
The parallel transport along a loop $\delta: [0,1]\arrow M$ such that $x=\delta(0)=\delta(1)$ gives an endomorphism $p_{\delta}\in\End(E_x)$ of the fiber $E_x$.  

\medskip
\begin{definition}
Let $\nabla: E\arrow\Lambda^1(M)\otimes E$ be a connection on a vector bundle $\pi: E\arrow G$.
    {\bf The holonomy group} of the connection $\nabla$  is the group of linear transformations of a fiber $E_x$ given by all parallel translations along all smooth loops based at $x$: $$\mathcal{H\!ol}_{\nabla}=\{p_{\delta}\,|\,p_{\delta}\in\End(E_x)\}.$$
\end{definition}

\medskip
\begin{definition} The holonomy group of a flat connection $\nabla$ is called {\bf the monodromy group of $\nabla$}. 
\end{definition}

\medskip
When $G$ is a Lie subgroup $\GL(n, K), K=\R$ or $\C$, its Lie algebra is canonically identified with the corresponding Lie subalgebra in $\Mat(n, K)$. In this case, recall that the exponential map $\exp:\g\arrow G$ is defined as follows
\begin{equation}\label{exp}
X\mapsto e^{tX}=\Id+\displaystyle\frac{t}{1!}X+\frac{t^2}{2!}X^2+\dots+\frac{t^k}{k!}X^k+\cdots,\,t\in\R, X\in\g,
\end{equation}
and its local inverse is given by a logarithm 
\begin{equation}\label{log}
   1+\gamma\mapsto \log(1+\gamma)=\gamma-\frac{\gamma^2}{2}+\frac{\gamma^3}{3}-\dots+(-1)^{k+1}\frac{\gamma^k}{k}+\dots,
\end{equation}
for all $\gamma\in U\subset G$ in a sufficiently small open neighborhood $U$ of the identity $e\in  G$. 

\medskip
Let $G$ be a simply connected nilpotent Lie group (it is diffeomorphic to a Euclidean
space $\R^n$ \cite[Corollary 2, p. 53]{VGO}) and $\Gamma$ a cocompact lattice in $G$, i.e. the left quotient $\Gamma\backslash G$ is a nilmanifold.

\medskip
\theorem\label{mongenlog}
Let $(\mathbb{B},\Grad)$ be a flat left equivariant vector bundle over a nilpotent Lie group $G$, $\Gamma$ a cocompact lattice and $(\mathbb{B}_{\Gamma},\Grad^{\Gamma})$ the corresponding flat bundle on the nilmanifold $\Gamma\backslash G$. Then the monodromy group $\mathcal{H\!ol}_{\Grad^{\Gamma}}$ of the connection $\Grad^{\Gamma}$ is generated by $e^{\log\gamma}$, where $\gamma\in\Gamma$. 

\proof
Let $\delta:[0,1]\arrow G$ be a $1$-parametric subgroup in a Lie group $G$, such that $\delta(0)=e$, $\delta(1)\in\Gamma$ and tangent to a vector $X=\log\gamma$. 
Then the solution of a linear differential equation $\nabla_{X}b=0$ along the curve $\gamma$ is given by the exponents $e^{tX}=e^{t\log\gamma}$. \endproof

\subsection{Maltsev completion}
In this section, we assume that a group $G$ is a finitely generated torsion-free nilpotent group.

\medskip
{\begin{definition}\label{Mal1}\cite{Mal}
 A nilpotent group $G$ is called {\bf Maltsev complete} if for each $g\in G$ and for all $n\in\Z_{>0}$ the equation $x^n=g$ has solutions in $G$. 
\end{definition}

\medskip
{\begin{definition}\cite{Mal}
 Let $\Gamma$ be a subgroup of a Maltsev complete nilpotent group ${G}$. Then the set $\hat{\Gamma}=\{g\in G\,|\, g^n\in\Gamma,\, n\in\Z\}\subset G$ is called {\bf the Maltsev completion} of a group $\Gamma$.
\end{definition}
 In \cite{Mal}, the definition of the Maltsev completion was given over the field $\Q$ but it could be done over any field $k$ of characteristic zero.

\medskip
\definition\label{Malcompl}\cite{GH}
Let $k$ be a field of characteristic zero. {\bf The Maltsev completion functor $\mathscr{M}_{k}$} is a functor from the category of finitely generated torsion-free nilpotent groups to the category of unipotent algebraic $k$-groups. If $G$ is a finitely generated nilpotent group, then $\mathscr{M}_{k}(G):=\overline{\Phi(G)}$, where $$\Phi:G\arrow\End(B)$$
is a faithful unipotent representation and $\overline{\Phi(G)}$ is the Zariski closure of the image, which lies in the subgroup of upper-triangular matrices.

\medskip
\remark
We will also denote the rational Maltsev completion of a group $G$ by $\hat{G}$, as it was introduced in \ref{Mal1}.

\medskip
Below we are listing some properties of the Maltsev completion \cite{GH}:

\medskip
\properties\label{prop}
\begin{enumerate}
\item[0.] A connected nilpotent Lie group is Maltsev complete;
\item[2.] Let $\Gamma$ be a subgroup of a nilpotent torsion-free Lie group $G$. By definition, $\mathscr{M}_{k}(\Gamma)\subset G$ is a Maltsev complete group. Its isomorphism class does not depend on $G$ and the embedding of $\Gamma$ to $G$;
\item[3.] The Maltsev completion $\mathscr{M}_{k}(\Gamma)$ is a minimal complete subgroup which contains $\Gamma$;
\item[4.] The Maltsev completion $\mathscr{M}_{k}(\Gamma)$ is equipped with a natural structure of an algebraic group over $k$.  
\end{enumerate}

\medskip
\claim\label{rem} 
The functor $\mathscr{M}_{\Q}$ 
provides a bijection between the finite-dimensional unipotent representations of $\Gamma$ over $\Q$ and the finite-dimensional $\Q$-representations of $\hat{\Gamma}$, considered as an algebraic group:
\[ 
\begin{tikzcd}
  \Gamma \arrow[r, hook, "\mathscr{M}_{\Q}"] \arrow[dr, "{\Phi}"]
    & \hat{\Gamma} \arrow[d, dashrightarrow, "\hat{\Phi}"]\\ &\End(B) \end{tikzcd}
\]
Moreover, the image $\Phi(\Gamma)$ is Zariski dense in $\hat{\Phi}(\hat{\Gamma})$.

\medskip
\proof \cite{GH}.
\endproof

\medskip
Note that in the context of \ref{mongenlog}, the holonomy group $\mathcal{H\!ol}_{\Grad^{\Gamma}}$ is isomorphic to the image of the representation of the lattice subgroup $\Gamma$:
\begin{equation}\label{gradrep}
\mathcal{H\!ol}_{\Grad^{\Gamma}} \cong \Phi(\Gamma).
\end{equation}

\medskip
The notion of a holonomy group has a geometric nature. We define an algebraic version of it as follows.

\medskip
\begin{definition}\label{alghol}
     Let $B$ be a vector space, $\g$ a Lie algebra, and $\nabla: B\arrow\g^*\otimes B$ an algebraic connection \eqref{nabla}. {\bf An algebraic holonomy group $ \mathcal{H\!ol}^a_\nabla$} is a subgroup of $\GL(B)$ generated by the matrix exponents:
\begin{equation*}
      \mathcal{H\!ol}^a_\nabla:=\langle e^{t\nabla_X}\,|\, t\in\R,\,\text{for all}\,X\in\g\,\rangle.
\end{equation*} 
\end{definition}

\medskip
\remark
We assume that an equivariant connection $\Grad$ on a $G$-equivariant vector bundle $B$ is flat. In this case, by \ref{repflat}
the associated algebraic connection defines a representation of a Lie algebra $X\mapsto\nabla_X\in\End(B)$.






\medskip
Let $G$ be a nilpotent group and $\hat{G}$ its rational Maltsev completion, which is a rational algebraic group, i.e. it can be written as $\hat{G}=\Spec(A)$,  where $A$ is its ring of regular functions. 

\medskip
\remark\label{realMalt} In \ref{Malcompl} we defined the functor of Maltsev completion $\mathscr{M}_{k}$ with coefficients in a field $k$.
{\bf The real Maltsev completion $\mathscr{M}_{\R}$} can be equivalently defined as follows: $$\hat{G}\otimes_{\Q}\R:=\Spec(A\otimes_{\Q}\R)=\mathscr{M}_{\R}(G).$$








\medskip
\theorem\label{ACL}
Let $(\mathbb{B},\Grad)$ be a flat left equivariant vector bundle over a nilpotent Lie group $G$, $\Gamma\subset G$ a cocompact lattice, and $(\mathbb{B}_{\Gamma},\Grad^{\Gamma})$ the induced bundle over the nilmanifold $\Gamma\backslash G$.
Suppose that $\nabla: B\arrow\g^*\otimes B$ is the algebraic connection associated with $\Grad$. Assume that the monodromy representation $\Phi: \Gamma\arrow\End(B)$ is unipotent.

Then
\begin{equation}
\overline{\mathcal{H\!ol}_{\Grad^{\Gamma}}}={\mathcal{H\!ol}}^a_{\nabla},
\end{equation}
where $\overline{\mathcal{H\!ol}_{\Grad^{\Gamma}}}\subseteq\GL(B)$ denotes the Zariski closure of the monodromy group of $(\mathbb{B}_{\Gamma},\Grad^{\Gamma})$.

\medskip
\proof 
Consider the chain of embeddings 
\begin{equation}
\mathcal{H\!ol}_{\Grad^{\Gamma}}\subseteq\widehat{\mathcal{H\!ol}_{\Grad^{\Gamma}}}\subseteq\widehat{\mathcal{H\!ol}_{\Grad^{\Gamma}}}\otimes_{\Q}\R, 
\end{equation}
where $\widehat{\mathcal{H\!ol}_{\Grad^{\Gamma}}}$ is the rational Maltsev completion of the group $\mathcal{H\!ol}_{\Grad^{\Gamma}}$.

Notice that $\mathcal{H\!ol}_{\Grad^{\Gamma}}\cong\Phi(\Gamma)$ is Zariski dense in $\widehat{\mathcal{H\!ol}_{\Grad^{\Gamma}}}$ by \cite[Lemma 3.11]{BrVer} and $\widehat{\mathcal{H\!ol}_{\Grad^{\Gamma}}}$ is Zariski dense in $\widehat{\mathcal{H\!ol}_{\Grad^{\Gamma}}}\otimes_{\Q}\R$.
Hence, 
\begin{equation}
\overline{\mathcal{H\!ol}_{\Grad^{\Gamma}}}=\widehat{\mathcal{H\!ol}_{\Grad^{\Gamma}}}\otimes_{\Q}\R.
\end{equation}

By \ref{alghol}, the algebraic holonomy group ${\mathcal{H\!ol}}^a_{\nabla}\subseteq\widehat{\mathcal{H\!ol}_{\Grad^{\Gamma}}}$ is a real group, which also contains $\mathcal{H\!ol}_{\Grad^{\Gamma}}$. Therefore, it has to be isomorphic to $\widehat{\mathcal{H\!ol}_{\Grad^{\Gamma}}}\otimes_{\Q}\R$, because of the minimality of the real Maltsev completion.

\endproof




\medskip
\corollary\label{AU} Let $\Gamma\backslash G$ be a hypercomplex nilmanifold with the flat Obata conection $\nabla^{Ob}$ \eqref{Obata} on $TG$. Then the action of the algebraic holonomy on $\g=\Lie(G)$ is unipotent.

\medskip
\proof
Let $\nabla$ be a flat algebraic connection on the Lie algebra $\g=T_eG$, associated with $\nabla^{Ob}$. Define the action $ \mathcal{H\!ol}^a_{\nabla}\times\g\arrow\g$ on $\g$ as follows:
\begin{equation}\label{action}
  (e^{\nabla_X},Y)\mapsto e^{\nabla_X}\cdot Y, 
\end{equation}
where $X,Y\in\g$.
By \ref{FGH}, the monodromy representation $\mathcal{H\!ol}_{\Grad^{\Gamma}}$ is unipotent. From \ref{ACL} it follows that $\mathcal{H\!ol}_{\Grad^{\Gamma}}$ is Zariski dense in ${\mathcal{H\!ol}}^a_{\nabla}$. Hence ${\mathcal{H\!ol}}^a_{\nabla}$ also unipotent.
\endproof

\subsection{Unipotent holonomy group and $\H$-solvability}

Let $(G, I, J, K)$ be a nilpotent Lie group with a left-invariant hypercomplex structure. In what follows we denote by $\nabla$ the Obata connection \eqref{Obata}, which we assumed to be flat, by $\mathcal{H\!ol}_\nabla$  and $\mathcal{H\!ol}^a_{\nabla}$ its holonomy (monodromy) and algebraic holonomy groups respectively.


\medskip
Let $\g$ be a nilpotent hypercomplex Lie algebra. Consider the smallest $\mathbb{H}$-invariant subspace of $\g$ containing the commutator ideal $[\g,\g]$:
\begin{equation*}
    \g_1^{\mathbb{H}}:=\H[\g,\g]=[\g,\g]+I[\g,\g]+J[\g,\g]+K[\g,\g].
\end{equation*}
Since $\H[\g,\g]$ contains the commutator ideal, it is an ideal, hence a subalgebra of $\g$.


\medskip
Define inductively the $\H$-invariant Lie subalgebras:
$$\g_i^{\H}:=\H[\g_{i-1}^{\H},\g_{i-1}^{\H}], i\in\Z_{\geq 2}.$$

\claim\label{preserve}
Let $(\g, I, J, K)$ be a hypercomplex nilpotent Lie algebra with the flat Obata connection. Then $\nabla_XY\in\g^{\H}_{i+1}$ for any $X,Y\in\g^{\H}_i$.

\medskip
\proof
Suppose that $X,Y\in\g^{\H}_i$. Then by \eqref{Obata} $\nabla_XY\in\g^{\H}_{i+1}$.
\endproof


\medskip
For each $i\in\Z_{\geq0}$, $X\in\g^{\H}_i$, the elements $e^{\nabla_X}\in\End({\g^{\H}_i})$ generate the algebraic holonomy group, associated to the Lie subalgebra $\g^{\H}_i$:
$$\mathcal{H\!ol}^a_{\nabla^i}=\langle e^{t\nabla_X}\,|\, t\in\R,\,X\in\g^{\H}_i\rangle.$$

\medskip
\remark\label{actionia}
Notice that $\mathcal{H\!ol}^a_{\nabla^i}$ acts on $\g^{\H}_i$ as in \eqref{action} and it could happen that it is not a subgroup or even a subset of the algebraic holonomy group $\mathcal{H\!ol}^a_{\nabla}$. Indeed, the 
action of $\nabla_X$ on $\g_i^{\H}$ could be trivial, but act non-trivially on $\g$ for some $X\in\g^{\H}_i$. However, the action of $\mathcal{H\!ol}^a_{\nabla^i}$ on 
$\g^{\H}_i$ coincides with the action of the corresponding subgroup of $\mathcal{H\!ol}^a_{\nabla}$. Therefore, it is unipotent on $\g^{\H}_i$ if $\mathcal{H\!ol}^a_{\nabla}$ is unipotent.


\medskip
The following theorem is the main result of this paper.

\medskip
\theorem
Let $(\g,I,J,K)$ be a hypercomplex nilpotent Lie algebra with the flat Obata connecton. Assume that the algebraic monodromy representation is unipotent. Then $\g$ is $\H$-solvable.

\medskip
\proof
The action of $\mathcal{H\!ol}^a_{\nabla^i}$ on the Lie subalgebra $\g^{\H}_i$ defined as in \eqref{action}. 
By the assumption, the representation of the algebraic holonomy on $\g$ is unipotent, hence, by \ref{actionia}, the action of $\mathcal{H\!ol}^a_{\nabla^i}$ is unipotent on the $\g^{\H}_i$.
 
Let  $\varphi:\mathcal{H\!ol}^a_{\nabla^i}\arrow\End((\g^{\H}_i)^*)$ be the action of the algebraic holonomy group on the dual Lie algebra $(\g_i^{\H})^*$. Every unipotent representation of a Lie algebra has an invariant vector. This implies the existence of a non-zero Obata-parallel 1-form for $\alpha\in\Lambda^1(\g_i^{\H})^*$. 
Again, we consider the holonomy action on $(\g_i^{\H})^*$ and for each $i$ we obtain a non-zero parallel (hence, closed\footnote{For any torsion-free connection $\nabla$ and any 1-form $\beta$ we have the equality $d\beta=\rm{Alt}(\nabla\beta)$, where $\rm{Alt}$ is the skew-symmetrization map.}) $1$-form $\alpha_i\in\Lambda^1(\g_i^{\H})^*$.

The intersection of kernels  $\Sigma_{\alpha_i}=\ker\alpha_i\cap\ker I\alpha_i\cap\ker J\alpha_i\cap\ker K\alpha_i$ gives an $\H$-invariant foliation, which contains $\g^{\H}_i$ as a proper subspace \eqref{intersectonofkernels}. 
This implies $\g_{i+1}^{\H}\subsetneq\g^{\H}_i$. The sequence $\g\subsetneq\g_1^{\H}\subsetneq\dots$ terminates in finitely many steps since $\g$ is finite-dimensional.
\endproof

\medskip
\corollary
Let $(N, I, J, K)$ be a hypercomplex nilmanifold with flat Obata connection $\nabla$. Then the Lie algebra $\g=\Lie(G)$ is $\H$-solvable.

\medskip
\proof
By \ref{AU}, the representation of the algebraic holonomy on $\g$ is unipotent.
\endproof


\noindent {\sc Yulia Gorginyan\\
{\sc Instituto Nacional de Matem\'atica Pura e
              Aplicada (IMPA) \\ Estrada Dona Castorina, 110\\
Jardim Bot\^anico, CEP 22460-320\\
Rio de Janeiro, RJ - Brasil\\
also:\\
{\sc Laboratory of Algebraic Geometry,\\
National Research University (HSE),\\
Department of Mathematics, 6 Usacheva Str.\\ Moscow, Russia}\\
\tt  iuliia.gorginian@impa.br }\\
}

\end{document}